\documentclass{amsart}

\usepackage[latin1]{inputenc}
\usepackage{subfigure}
\usepackage[english]{babel}
\usepackage{amsmath}
\usepackage{amssymb}
\usepackage{centernot}
\usepackage{etoolbox}
\usepackage{titlesec}
\patchcmd{\section}{\scshape}{\bfseries}{}{}
\makeatletter
\renewcommand{\@secnumfont}{\bfseries}
\makeatother
\newcommand*{\justifyheading}{\raggedright}
\titleformat {\section}
  {\normalfont \Large \bfseries \justifyheading}{\thesection.}{1em}{}
\titleformat {\subsection}
  {\normalfont \large \bfseries \justifyheading}{\thesubsection.}{1em}{}
\patchcmd{\abstract}{\scshape\abstractname}{\textbf{\abstractname}}{}{}
\numberwithin{equation}{section}
\newtheorem{thm}{Theorem}[section]
\newtheorem{prop}[thm]{Proposition}
\newtheorem{lem}[thm]{Lemma}
\newtheorem{cor}[thm]{Corollary}

\newtheorem{obs}[thm]{Observation}
\newtheorem{theo}[thm]{Theorem}
\newtheorem{defn}[thm]{Definition}

\newtheorem{conv}[thm]{Convention}
\theoremstyle{definition}

\newtheorem{prob}[thm]{Problem}
\newtheorem{exam}[thm]{Example}
\newenvironment{proof2}[1][Proof.]{\begin{trivlist}
\item[\hskip \labelsep {\bfseries \itshape #1}]}{\end{trivlist}}
\newenvironment{proof3}[1][Proof.]{\begin{trivlist}
\item[\hskip \labelsep {\itshape #1}]}{\end{trivlist}}
\newcommand{\bij}[2]{\text{Bij$\left(#1,#2\right)$}}
\newcommand{\rev}{\text{rev}}
\newcommand{\expe}[1]{\text{\bf{E$\left[\right.$}}#1\text{\bf{$\left.\right]$}}}
\newcommand{\proba}[1]{\text{\bf{Pr$\left[\right.$}}#1\text{\bf{$\left.\right]$}}}
\newcommand{\Ou}{O_{\text{up}}}
\newcommand{\Od}{O_{\text{down}}}
\newcommand{\Img}[1]{\text{Im}\left(#1\right)}
\newcommand{\Chn}[1]{\mathcal{C}\left(#1\right)}
\newcommand{\Ord}[1]{\mathcal{O}\left(#1\right)}
\newcommand{\tri}[1]{\emph{tri}\left(#1\right)}
\newcommand{\inc}[1]{\emph{incom}\left(#1\right)}
\newcommand{\com}[1]{\emph{com}\left(#1\right)}
\newcommand{\R}{\mathbb{R}}
\newcommand{\N}{\mathbb{R}_{\geq 0}}

\newcommand{\indeg}[1]{\text{indeg}\left(#1\right)}
\newcommand{\outdeg}[1]{\text{outdeg}\left(#1\right)}
\newcommand{\Arr}{\mathcal{A}}
\newcommand{\Vol}[1]{\text{Vol}\left(#1\right)}
\newcommand{\Voli}[1]{\text{Vol}^\circ\left(#1\right)}
\newcommand{\stab}[1]{\text{STAB}\left(#1\right)}
\newcommand{\card}[1]{\left|#1\right|}

\newcommand{\set}[1]{\left\{#1\right\}}
\newcommand{\Frac}[2]{\displaystyle\frac{#1}{#2}}
\newcommand{\rgr}[2]{G_{#1,#2}}
\newcommand{\logt}{\log_2}

\usepackage{graphicx}
\usepackage{mathrsfs}
\usepackage[foot]{amsaddr}

\usepackage[round]{natbib}

\title[Graph Orientations and Linear Extensions.]{Graph Orientations and
Linear Extensions.}
\author{Benjamin Iriarte}
\address{Department of Mathematics, Massachusetts Institute of Technology, Cambridge MA, 02139, USA}
\email{biriarte@math.mit.edu}

\keywords{acyclic orientation, linear extension, poset, comparability
graph, stable polytope.}

\begin{document}
\begin{abstract}
Given an underlying undirected simple graph, we consider the set of all
acyclic orientations of its edges. Each of these orientations induces a
partial order on the vertices of our graph and, therefore, we can count the
number of linear extensions of these posets. We want to know which choice
of orientation maximizes the number of linear extensions of the
corresponding poset, and this problem will be solved essentially for
comparability graphs and odd cycles, presenting several proofs. 
The corresponding enumeration problem for
arbitrary simple graphs will be studied, including the case of random graphs; 
this will culminate 
in 1) new bounds for the volume of the stable polytope and 2) strong concentration
results for our main statistic and for the graph
entropy, which hold true $a.s.$ for random graphs. We will then argue that our problem
springs up naturally in the theory of graphical arrangements and graphical zonotopes. 
\end{abstract}

\maketitle

\section{Introduction.}

Linear extensions of partially ordered sets have been the object of much
attention and their uses and applications remain increasing. Their number
is a fundamental statistic of posets, and they relate to ever-recurring
problems in computer science due to their role in sorting problems. Still,
many fundamental questions about linear extensions are unsolved, including
the well-known \emph{$1\slash 3$-$2\slash 3$ Conjecture}. 
Efficiently enumerating linear extensions of certain posets is difficult, 
and the general problem has been 
found to be $\sharp$P-complete in \cite{nph}. 

Directed acyclic graphs, and similarly, acyclic orientations of simple
undirected graphs, are closely related to posets, and their
problem-modeling values in several disciplines, including the biological
sciences, 
needs no introduction. We
propose the following problem:
\begin{prob}
Suppose that there are $n$ individuals with a known contagious disease,
and suppose that we know which pairs of these individuals were in the same
location at the same time. Assume that at some initial points, some of the
individuals fell ill, and then they started infecting other people and so
forth, spreading the disease until all $n$ of them were infected. Then,
assuming no other knowledge of the situation, what is the most likely way
in which the disease spread out?        
\end{prob}

Suppose that we have an underlying connected undirected simple graph
$G=G(V,E)$ with $n$ vertices. If we first pick uniformly at random a
bijection $f:V\rightarrow [n]$, and then orient the edges of $E$ so that
for every $\{u,v\}\in E$ we select $(u,v)$ (read $u$ directed to $v$)
whenever $f(u)<f(v)$, we obtain an acyclic orientation of $E$. In turn,
each acyclic orientation induces a partial order on $V$ in which $u<v$ if
and only if there is a directed path $(u,u_1),(u_1,u_2),\dots,(u_k,v)$ in
the orientation. In general, several choices of $f$ above will
result in the same acyclic orientation. However, the most likely acyclic
orientations so obtained will be the ones whose induced posets have the
maximal number of linear extensions, among all posets arising from acyclic
orientations of $E$. Our problem then becomes that of deciding which
acyclic orientations of $E$ attain this optimality property of maximizing
the number of linear extensions of induced posets. This problem, referred to throughout this
article as the \emph{main problem for $G$}, was
raised by \cite{ksaito} for the case of trees, yet, a
solution for the case of bipartite graphs had been obtained already by
\cite{stach}. The main problem brings up the natural associated 
enumerative question: For a graph $G$, what is the
maximal number of linear extensions of a poset induced by an acyclic orientation of $G$? This statistic 
for simple graphs will be herein referred to as the \emph{main statistic} (Definition~\ref{def:maxi}).  

The central goal of this initial article on the subject will be to begin a rigorous study of the main problem 
from the points of view of structural and enumerative combinatorics. We will introduce 
{\bf 1)} techniques to find optimal orientations of graphs that are provably correct for certain families of graphs, 
and {\bf 2)} techniques to estimate the main statistic for more general 
classes of graphs and to further understand aspects of its distribution across all graphs.   

In Section~\ref{sec:intro}, we will present an
elementary approach to the main problem for both bipartite graphs and odd cycles. 
This will serve as motivation and preamble for the remaining sections. 
In particular, in Section~\ref{sec:bip}, a new 
solution to the main problem for bipartite graphs will be obtained, different to that of \cite{stach} in
that we explicitly construct a function
that maps injectively linear extensions of non-optimal acyclic
orientations to linear extensions of an optimal orientation. As we will observe, optimal orientations
of bipartite graphs are precisely the \emph{bipartite orientations} (Definition~\ref{def:bipor}). 
Then, in Section~\ref{sec:oc}, we will
extend our solution for bipartite graphs to odd cycles, proving that optimal orientations
of odd cycles are precisely the \emph{almost bipartite orientations} (Definition~\ref{def:abip}). 

In Section~\ref{sec:cg}, we will introduce two new
techniques, one geometrical and the other
poset-theoretical, that lead to different solutions for the case of
comparability graphs. 
Optimal orientations of comparability graphs
are precisely the \emph{transitive orientations} (Definition~\ref{def:trans}), a result that
generalizes the solution for bipartite graphs. 
The techniques developed on Section~\ref{sec:cg} will allow
us to re-discover the solution for odd cycles and to state 
inequalities for the general enumeration
problem in Section~\ref{sec:fe}. The recurrences for the 
number of linear extensions of posets presented in
Corollary~\ref{cor:tri} had been
previously established in \cite{ehs} using {\it promotion and evacuation}
theory, but we will obtain them independently
as by-products of certain network flows in Hasse diagrams. Notably, \cite{stach}
had used some instances of these recurrences to solve the main problem for
bipartite graphs. 

Further on, in Section~\ref{sec:fe}, we will
also consider the enumeration problem for the case of random graphs
with distribution $\rgr{n}{p},0<p<1$, and obtain 
tight concentration results for our main statistic, across all graphs. Incidentally, this will lead to new inequalities 
for the volume of the stable polytope and to a very strong concentration result for the graph entropy
(as defined in \cite{lov}), which hold $a.s.$ for random graphs.  

Lastly, in Section~\ref{sec:ft}, we will
show that the main problem for a graph arises naturally from the corresponding graphical arrangement by asking
for the regions with maximal \emph{fractional volume} (Proposition~\ref{prop:vols}). More surprisingly, we will also observe that the solutions to the main problem
for comparability graphs and odd cycles correspond to certain vertices of the corresponding graphical zonotopes (Theorem~\ref{theo:cgz}). 

\begin{conv}
Let $G=G(V,E)$ be a simple undirected graph. Formally, an \emph{orientation $O$ of $E$ (or $G$)} 
is a map $O:E\rightarrow V^2$ such that for all $e:=\{u,v\}\in E$, we have 
$O(e)\in\{(u,v),(v,u)\}$. Furthermore, $O$ is said to be \emph{acyclic} if 
the directed graph on vertex set $V$ and directed-edge set $O(E)$ is acyclic. 
On numerous occasions, we will somewhat abusively also identify
an acyclic orientation $O$ of $E$ with the set $O(E)$, or with the poset that
it induces on $V$, doing this with the aim to reduce extensive wording.   
\end{conv}

When defining posets herein, we will also try to make clear the
distinction between the ground set of the poset and its order relations.

\section{Introductory results.}\label{sec:intro}
\subsection{The case of bipartite graphs.}\label{sec:bip}

The goal of this section is to present a combinatorial proof that the
number of linear extensions of a bipartite graph $G$ is maximized when we
choose a \emph{bipartite orientation} for $G$. Our method is to find an
injective function from the set of linear extensions of any fixed acyclic
orientation to the set of linear extensions of a bipartite
orientation, and then to show that this function is not surjective whenever
the initial orientation is not bipartite. Throughout the section, let $G$
be bipartite with $n\geq 1$ vertices. 

\begin{defn}\label{def:bipor}
Suppose that $G=G(V,E)$ has a bipartition $V=V_1\sqcup V_2$. Then, the
orientations that either choose $(v_1,v_2)$ for all $\{v_1,v_2\}\in E$ with
$v_1\in V_1$ and $v_2\in V_2$, or
$(v_2,v_1)$ for all $\{v_1,v_2\}\in E$ with $v_1\in V_1$ and $v_2\in V_2$,
are called \emph{bipartite orientations} of $G$. 
\end{defn}
\begin{defn}
For a graph $G$ on vertex set $V$ with $|V|=n$, we will denote by
$\bij{V}{[n]}$
the set of bijections from $V$ to $[n]$.
\end{defn}

As a training example, we consider the case when we transform linear
extensions of one of the bipartite orientations into linear extensions of
the other bipartite orientation. We expect to obtain a bijection for this
case.

\begin{prop}\label{prop:bip}
Let $G=G(V,E)$ be a simple connected undirected bipartite graph, with
$n=|V|$. 
Let $\Od$ and $\Ou$ be the two bipartite orientations of $G$. Then, there
exists a bijection between the set of linear extensions of 
$\Od$ and the set of linear extensions of $\Ou$.
\end{prop}
\begin{proof2}
Consider the automorphism $\rev$ of the set $\bij{V}{[n]}$ 
given by $\rev(f)(v)=n+1-f(v)$ for all $v\in V$ and $f\in\bij{V}{[n]}$.
It is clear that $(\rev\circ\rev)(f)=f$. However, since
$f(u)>f(v)$ implies $\rev(f)(u)<\rev(f)(v)$, then $\rev$ reverses all
directed paths in any $f$-induced acyclic orientation of $G$, and in
particular the restriction of $\rev$
to the set of linear extensions of $\Od$ has image $\Ou$, and viceversa.  
\par \qed \end{proof2}

We now proceed to study the case of general acyclic orientations of the
edges of $G$. 
Even though similar in flavour to Proposition~\ref{prop:bip}, our new
function will not in general correspond to the function presented in the
proposition when restricted to the case of bipartite orientations. 

To begin, we define the main automorphisms of $\bij{V}{[n]}$ that will
serve as building blocks for constructing the new function. 

\begin{defn}\label{def:auto}
Consider a simple graph $G=G(V,E)$ with $|V|=n$. 
For different vertices $u,v\in V$, let $\rev_{uv}$ be the automorphism of
$\bij{V}{[n]}$ given by the following rule: For all $f\in \bij{V}{[n]}$,
let
$$\begin{array}{l l}
\rev_{uv}(f)(u) & =f(v),\\
\rev_{uv}(f)(v) & =f(u),\\
\rev_{uv}(f)(w) & =f(w) \text{ if $w\in V\backslash\{u,v\}$.}
  \end{array}$$
\end{defn}

It is clear that $(\rev_{uv}\circ\rev_{uv})(f)=f$ for all
$f\in\bij{V}{[n]}$. Moreover, we will need the following technical
observation about $\rev_{uv}$.   

\begin{obs}\label{obs:tec}
Let $G=G(V,E)$ be a simple graph with $|V|=n$ and consider a bijection 
$f\in\bij{V}{[n]}$. 
Then, if for some $u,v,x,y\in V$ with $f(u)<f(v)$ we have 
that $\rev_{uv}(f)(x)>\rev_{uv}(f)(y)$ but $f(x)<f(y)$, then
$f(u)\leq f(x)<f(y)\leq f(v)$ and furthermore,
at least one of $f(x)$ or $f(y)$ must be equal to one of $f(u)$ or $f(v)$.
\end{obs}

Let us present the main result of this section, obtained based on the
interplay between acyclic orientations and bijections in $\bij{V}{[n]}$.

\begin{theo}\label{thm:bip}
Let $G=G(V,E)$ be a connected bipartite simple graph with $|V|=n$, and
with bipartite orientations $\Od$ and $\Ou$. Let also $O$ be an acyclic
orientation of $G$. Then, there exists an injective function $\Theta$ from
the set of linear extensions of $O$ to the set of linear extensions of
$\Ou$ and furthermore, $\Theta$ is surjective if and only if $O= \Ou$ or
$O= \Od$.  
\end{theo}
\begin{proof2}
Let $f$ be a linear extension of $O$, and without loss of
generality assume that $O\neq \Ou$. We seek to find a function $\Theta$
that transforms $f$ into a linear extension of $\Ou$ injectively. 
The idea will be to describe how $\Theta$ acts on $f$ as a composition of
automorphisms of the kind presented in Definition~\ref{def:auto}. Now, we
will find the terms of the composition in an inductive way, where at each
step we consider the underlying configuration obtained after the previous
steps. In particular, the choice of terms in the composition will depend on
$f$.  The inductive steps will be indexed using a positive integer variable $k$,
starting from $k=1$, and at each step we will know an acyclic orientation
$O_k$ of $G$, a set $B_k$ and a function $f_k$. The set $B_k\subseteq V$
will always be defined as the set of all vertices incident to an edge whose
orientation in $O_k$ and $\Ou$ differs, and $f_k$ will be a particular
linear extension of $O_k$ that we will define.  

Initially, we set $O_1=O$ and $f_1=f$, and calculate $B_1$. Now, suppose
that for some fixed $k\geq 1$ we know $O_k, B_k$ and $f_k$, and we want to
compute $O_{k+1}, B_{k+1}$ and $f_{k+1}$. If $B_k=\emptyset$, then
$O_k=\Ou$ and
$f_k$ is a linear extension of $\Ou$, so we stop our recursive process. 
If not, then $B_k$ contains elements $u_k$ and $v_k$ such that $f_k(u_k)$
and $f_k(v_k)$ are respectively minimum and maximum elements of
$f_k(B_k)\subseteq[n]$. Moreover, $u_k\neq v_k$. We will then let
$f_{k+1}:=\rev_{u_kv_k}(f_k)$, $O_{k+1}$ be the acyclic orientation of $G$
induced by $f_{k+1}$, and calculate $B_{k+1}$ from $O_{k+1}$.

If we let $m$ be the minimal positive integer for which
$B_{m+1}=\emptyset$, then
$\Theta(f)=(\rev_{u_mv_m}\circ\cdots\circ\rev_{u_2v_2}\circ\rev_{u_1v_1})(f)$.
The existence of $m$ follows from observing that $B_{k+1}\subsetneq B_{k}$
whenever $B_k\neq\emptyset$. In particular, if $B_k\neq\emptyset$, then
$u_k,v_k\in B_{k}\backslash B_{k+1}$ and so $1\leq m\leq \left\lfloor
\frac{|B_1|}{2}\right\rfloor$. It follows that the pairs
$\left\{\{u_k,v_k\}\right\}_{k\in[m]}$ are pairwise disjoint,
$f(u_k)=f_k(u_k)$ and
$f(v_k)=f_k(v_k)$ for all $k\in[m]$, and
$f(u_1)<f(u_2)<\dots<f(u_m)<f(v_m)<\dots<f(v_2)<f(v_1)$. As a consequence,
the 
automorphisms in the composition description of $\Theta$ commute. Lastly,
$f_{m+1}$ will
be a linear extension of $\Ou$ and we stop the inductive process by
defining
$\Theta(f)=f_{m+1}$. 

To prove that $\Theta$ is injective, note that given $O$ and $f_{m+1}$ as
above, we can recover uniquely $f$ by imitating our procedure to find
$\Theta(f)$. Firstly, set $g_1:=f_{m+1}$ and $Q_1:=\Ou$, and compute
$C_1\subseteq V$ as the set of vertices incident to an edge whose
orientation differs in $Q_1$ and $O$. Assuming prior knowledge of
$Q_k, C_k$ and $g_k$, and whenever $C_k\neq\emptyset$ for some positive
integer $k$, find the elements of $C_k$ whose images under $g_k$ are
maximal and minimal in $g_k(C_k)$. By the discussion above and
Observation~\ref{obs:tec}, we check that these are respectively and
precisely $u_k$ and $v_k$. Resembling the previous case, we will then let
$g_{k+1}:=\rev_{u_kv_k}(g_k)$, $Q_{k+1}$ be the acyclic orienation of $G$
induced by $g_k$, and compute
$C_{k+1}$ accordingly as the set of vertices incident to an edge with
different orientation in $Q_{k+1}$ and $O$. Clearly $g_{m+1}=f$, and the
procedure shows that $\Theta$ is invertible in its image. 

To establish that $\Theta$ is not surjective whenever $O\neq \Od$, note
that then $O$ contains a directed $2$-path $(w,u)$ and $(u,v)$. Without
loss of generality, we may assume that the orientation of these edges in
$\Ou$ is given by 
$(w,u)$ and $(v,u)$. But then, a linear extension $g$ of $\Ou$ in which
$g(u)=n$ and $g(v)=1$ is not in $\Img{\Theta}$ since otherwise, 
using the notation and framework discussed above, there would exist 
different $i,j\in[m]$ such that $u_i=u$ and $v_j=v$, which then
contradicts the choice of $u_1$ and $v_1$. This completes the proof.   
\par \qed \end{proof2}

\begin{figure}
\centering
\includegraphics[width=0.7\textwidth]{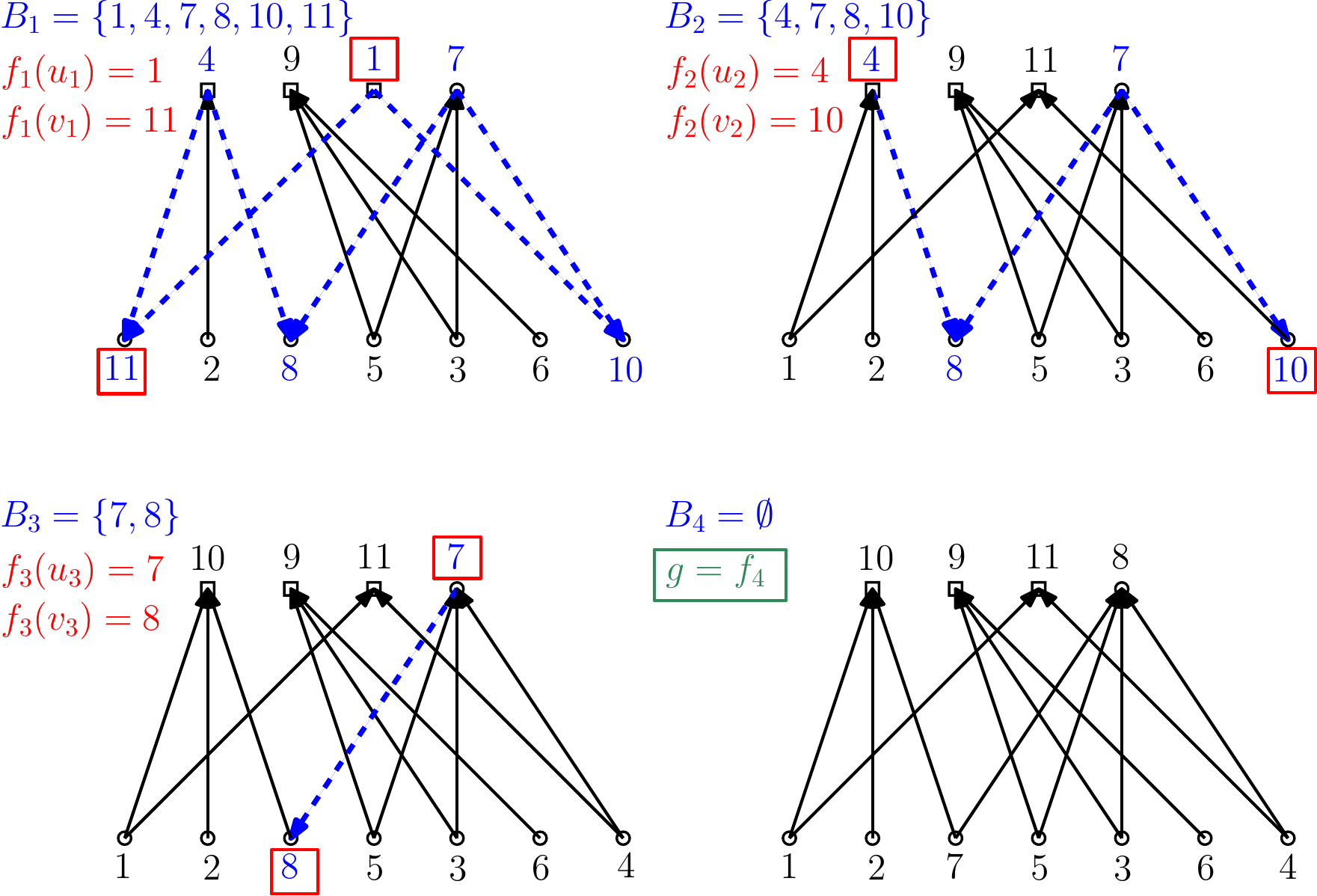}
\caption{\label{fig:bip}An example of the function $\Theta$ for the case
of bipartite graphs. Squares
show the numbers that will be flipped at each step, and dashed arrows
indicate arrows whose orientation
still needs to be reversed.}
\end{figure}

\subsection{Odd cycles.}\label{sec:oc}

In this section $G=G(V,E)$ will be a cycle on $2n+1$ vertices with $n\geq
1$. 
The case of odd cycles follows as an immediate extension of the case of
bipartite graphs, but it will also be covered under a different guise in
Section~\ref{sec:fe}. As expected, the acyclic orientations of the edges of
odd cycles that maximize the number of linear extensions resemble as much
as possible bipartite orientations. This is now made precise.

\begin{defn}\label{def:abip}
For an odd cycle $G=G(V,E)$, we say that an ayclic orientation of its
edges is \emph{almost bipartite} if under the orientation there exists exactly one
directed 
$2$-path, i.e. only one instance of $(u,v)$ and $(v,w)$ with $u,v,w\in V$.
\end{defn}

\begin{theo}\label{th:oc}
Let $G=G(V,E)$ be an odd cycle on $2n+1$ vertices with $n\geq 1$. Then,
the acyclic orientations of $E$ that maximize the number of linear
extensions are the almost bipartite orientations.
\end{theo}
\begin{proof2}[First proof.]
Since the case when $n=1$ is straightforward let us assume that $n\geq 2$,
and consider an arbitrary acyclic orientation $O$ of $G$.
Again, our method will be to construct an injective function $\Theta'$
that transforms every linear extension of $O$ into a linear extension of
some fixed almost bipartite orientation of $G$, where the specific choice
of almost bipartite orientation will not matter by the symmetry of $G$.

To begin, note that there must exist a directed $2$-path in $O$, say
$(u,v)$ and $(v,w)$ for some $u,v,w\in V$. Our goal
will be to construct $\Theta'$ so that it maps into the set of linear
extensions of the almost bipartite orientation $O_{uvw}$ in which our
directed path $(u,v),(v,w)$ is the unique directed $2$-path.
To find $\Theta'$, first consider the bipartite graph
$G'$ with vertex set $V\backslash\{v\}$ and edge set
$E\backslash\left(\{u,v\}\cup\{v,w\}\right)\cup\{u,w\}$, along with 
the orientation $O'$ of its edges that agrees on common edges 
with $O$ and contains $(u,w)$. Clearly $O'$ is acyclic. If $f$ is a linear
extension of $O$, we regard the restriction $f'$ of $f$ to
$V\backslash\{v\}$ as a strict order-preserving map on $O'$, and
analogously to the proof of Theorem~\ref{thm:bip}, we can transform
injectively $f'$ into a strict order-preseving map $g'$ with
$\Img{g'}=\Img{f'}=\Img{f}\backslash\{f(v)\}$ of the bipartite orientation
of $G'$ that contains $(u,w)$. Now, if we define $g\in\bij{V}{[n]}$ via
$g(x)=g'(x)$ for all $x\in V\backslash\{v\}$ and $g(v)=f(v)$, we see that
$g$ is a linear extension of $O_{uvw}$. We let $\Theta'(f)=g$. 

The technical work for proving the general injectiveness of $\Theta'$, and
its non-surjectiveness when $O$
is not almost bipartite, has already been presented in the proof of 
Theorem~\ref{thm:bip}: That $\Theta'$ is injective follows from the
injectiveness of the map transforming $f'$ into $g'$, and then by noticing
that $f(v)=g(v)$. Non-surjectiveness follows from noting that if $O$ is not
almost bipartite, then $O$ contains a directed $2$-path $(a,b),(b,c)$ with
$a,b,c\in V$ and $b\neq v$, so we cannot have simultaneously $g'(a)=
\min{\Img{f'}}$ and $g'(c)=\max{\Img{f'}}$.   
\par \qed \end{proof2}

\begin{figure}
\centering
\includegraphics[width=0.7\textwidth]{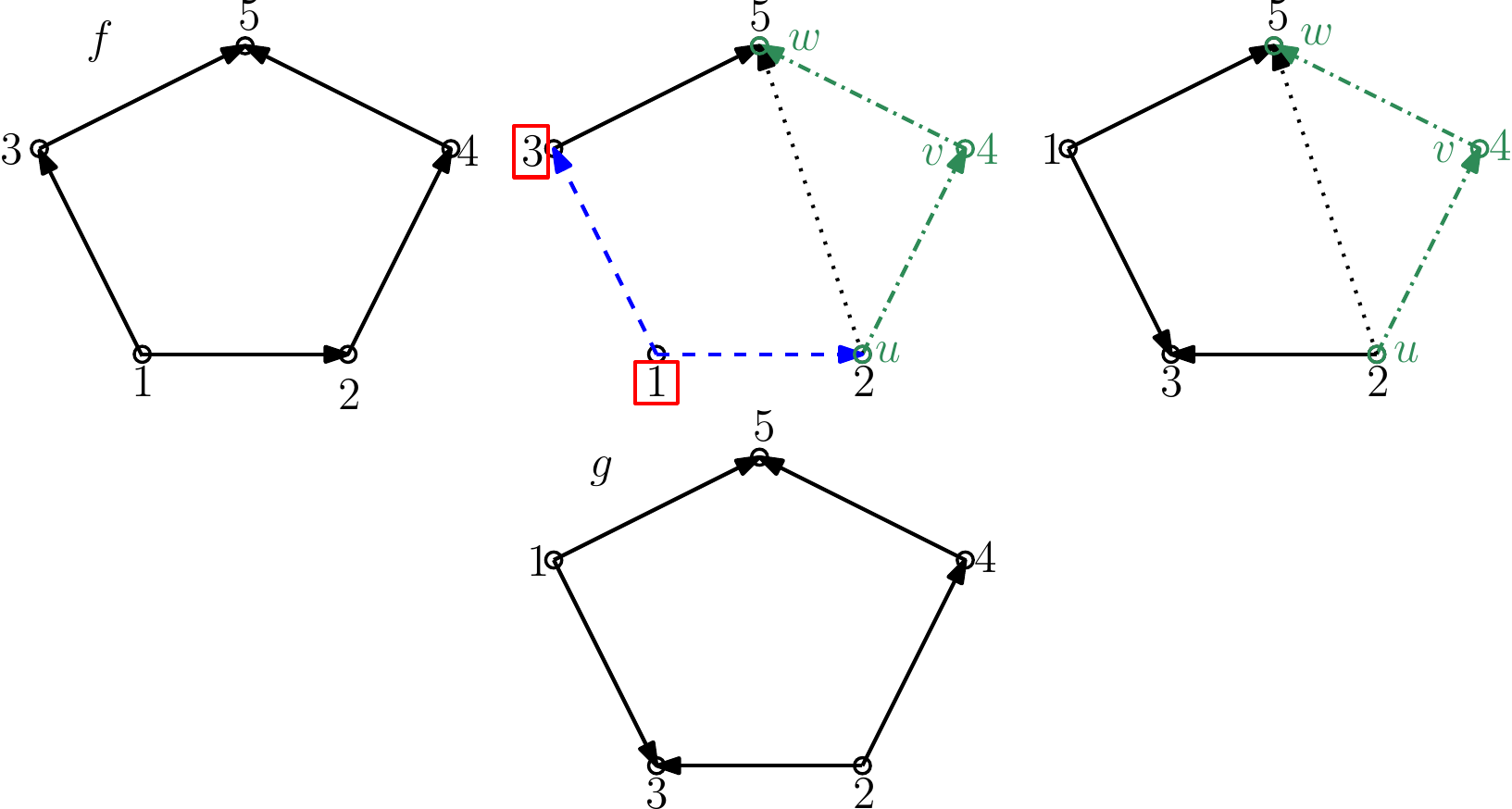}
\caption{\label{fig:odd}An example of the function $\Theta'$ for the case
of odd cycles. Squares show the numbers 
that will be flipped at each step. Dashed arrows indicate arrows whose
orientation still needs to be reversed, while
dashed-dotted arrows indicate those whose orientation will never be reversed. In
particular, $4$ will remain labeling the same
vertex during all steps.}
\end{figure}

\section{Comparability graphs.}\label{sec:cg}

In this section, we will study our main problem using more general
tehniques. As a consequence, we will be able to understand the case
of comparability graphs, which includes bipartite graphs as a special
case. Let us first recall the main object of this section:
\begin{defn}\label{def:compar}
A \emph{comparability graph}
is a simple undirected graph $G=G(V,E)$ for which there exists a partial
order on $V$ under which
two different vertices $u,v\in V$ are comparable if and only if
$\{u,v\}\in E$.  
\end{defn}

The acyclic orientations of the edges of a comparability graph $G$ that
maximize the number of linear extensions are precisely the orientations
that induce posets whose comparability graph agrees with $G$.  

Comparability graphs have been largely
discussed in the literature, mainly due to their connection with
partial orders and because they are perfectly orderable graphs and more
generally, perfect graphs. Comparability graphs, perfectly orderable
graphs
and perfect graphs are all large hereditary classes of graphs. In Gallai's
fundamental work in  \cite{pg}, a characterization of comparability graphs
in terms of
forbidden subgraphs was given and the concept of {\it modular
decomposition} of a graph was introduced.

Note that, given a comparability graph $G=G(V,E)$, we can find at least two
partial
orders on $V$ induced by acyclic orientations of $E$ whose comparability
graphs (obtained as discussed above) agree precisely with $G$, and the
number of such posets depends on the modular structure of $G$. Let us
record this idea in a definition.

\begin{defn}\label{def:trans}
Let $G=G(V,E)$ be a comparability graph, and let $O$ be an acyclic
orientation of $E$ such that the comparability graph of the partial order
of $V$ induced by $O$ agrees precisely with $G$. Then, we will say that $O$
is a \emph{transitive orientation} of $G$.
\end{defn}

We will present two methods for proving our main result. The first one (Subsection~\ref{subsec:g})
relies on Stanley's transfer map between the \emph{order polytope} and the
\emph{chain polytope} of a poset, and the second one (Subsection~\ref{subsec:pt}) is made possible by
relating our problem to \emph{network flows}. 

\subsection{Geometry.}\label{subsec:g}

To begin, let us recall the main definitions and notation related to the
first method. 

\begin{defn}
We will consider $\R^n$ with euclidean topology, and let
$\{e_j\}_{j\in[n]}$ be the standard basis of $\R^n$. For $J\subseteq[n]$,
we will define $e_{J}:=\sum_{j\in J}e_j$ and $e_{\emptyset}:=0$;
furthermore, for
$x\in\R^n$ we will let $x_J:=\sum_{j\in J}x_j$ and $x_{\emptyset}:=0$.
\end{defn}
\begin{defn}\label{def:cop}
Given a partial order $P$ on $[n]$, the \emph{order polytope} of $P$ is
defined as: $$\Ord{P}:=\left\{ x\in\R^n:0\leq x_i\leq 1\text{ and }x_j\leq
x_k\text{ whenever }j\le_Pk\text{, }\forall\text{ }i,j,k\in [n]\right\}.$$
The \emph{chain polytope} of $P$ is defined as:
$$\Chn{P}:=\{x\in\R^n:x_i\geq 0\text{, $\forall$ $i\in [n]$ and }x_C\leq
1\text{ whenever $C$ is a chain in $P$}\}.$$
\emph{Stanley's transfer map} $\phi:\Ord{P}\rightarrow\Chn{P}$ is the function
given by: 
$$\phi(x)_i=\left\{\begin{array}{cc}
x_i-\max_{j\lessdot_P i}x_j&\text{ if $i$ is not minimal in $P$,}\\
x_i&\text{ if $i$ is minimal in $P$.}
\end{array}\right.$$
\end{defn}

Let $P$ be a partial order on $[n]$.
It is easy to see from the definitions that the vertices of $\Ord{P}$ are
given 
by all the $e_I$ with $I$ an order filter of $P$, and those of $\Chn{P}$
are given by all the $e_A$ with $A$ an antichain of $P$.

Now, a well-known result of~\cite{stp} states that
$\Vol{\Ord{P}}=\frac{1}{n!}e(P)$ where $e(P)$ is the number of linear
extensions of $P$. This result can be proved by considering the unimodular
triangulation of $\Ord{P}$ whose maximal (closed) simplices have the form
$\Delta_{\sigma}:=\{x\in\R^n:0\leq x_{\sigma^{-1}(1)}\leq
x_{\sigma^{-1}(2)}\leq\dots\leq 
x_{\sigma^{-1}(n)}\leq 1\}$ with $\sigma:P\rightarrow\mathbf{n}$ a linear
extension of $P$. However, the volume of 
$\Chn{P}$ is not so direct to compute. To find $\Vol{\Chn{P}}$ Stanley
made use of the transfer map $\phi$, a pivotal idea that we now wish to
describe in detail since it will provide a geometrical point of view on our
main problem.

It is easy to see that $\phi$ is invertible and its inverse
can be described by: 
$$\phi^{-1}(x)_i=\max_{\substack{C\text{ chain in $P$:}\\ 
i\text{ is maximal in $C$}}}x_C, \text{ for all $i\in[n]$ and
$x\in\Chn{P}$}.$$
As a consequence, we see that $\phi^{-1}(e_A)=e_{A^{\vee}}$ for all
antichains
$A$ of $P$, where $A^{\vee}$ is the order filter of $P$ induced by $A$.
It is also straightforward to notice that $\phi$ is linear on each of the
$\Delta_{\sigma}$ with $\sigma$ a linear extension of $P$, by staring at
the definition of $\Delta_{\sigma}$. Hence, for fixed $\sigma$ and for each
$i\in[n]$, we can consider the
order filters $A_{i}^{\vee}:=\sigma^{-1}([i,n])$ along with their
respective minimal elements 
$A_i$ in $P$, and notice that $\phi(e_{A_{i}^{\vee}})=e_{A_i}$ and also
that 
$\phi(0)=0$. From there, $\phi$ is now easily 
seen to be a unimodular linear map on
$\Delta_{\sigma}$, and so
$\Vol{\phi\left(\Delta_{\sigma}\right)}=\Vol{\Delta_{\sigma}}=\frac{1}{n!}$.
Since $\phi$ is invertible, without unreasonable effort we have obtained
the following central result:

\begin{theo}[\cite{stp}]\label{thm:stanley}
Let $P$ be a partial order on $[n]$. Then, 
$\Vol{\Ord{P}}=\Vol{\Chn{P}}=\frac{1}{n!}e(P)$, where $e(P)$ is the number
of linear extensions of $P$.
\end{theo}

\begin{defn}\label{defn:stabpol}
Given a simple undirected graph $G=G([n],E)$, the
\emph{stable polytope} $\stab{G}$ of $G$ is the full dimensional polytope
in $\R^n$ obtained as the convex hull of all the vectors $e_I$,
where $I$ is a stable (a.k.a. independent) set of $G$. 
\end{defn}

Now, the chain polytope of a partial order $P$ on $[n]$ is clearly the
same as the stable polytope $\stab{G}$ of its comparability graph
$G=G([n],E)$ since antichains of $P$ correspond to stable sets of $G$. In
combination with Theorem~\ref{thm:stanley}, this shows that the number of
linear extensions is a comparability invariant, i.e. two posets with
isomorphic comparability graphs have the same number of linear extensions. 

We are now ready to present the first proof of the main result for
comparability graphs. We will assume connectedness of $G$ for
convenience in the presentation of the second proof.
\begin{theo}\label{th:cg}
Let $G=G(V,E)$ be a connected comparability graph. Then, the acyclic
orientations of $E$ that maximize the number of linear extensions are
exactly the transitive orientations of $G$.    
\end{theo}
\begin{proof2}[First proof.]
Without loss of generality, assume that $V=[n]$. Let $O$ be an acyclic
orientation of $G$ inducing a partial order $P$ on $[n]$. 
If two vertices $i,j\in[n]$ are incomparable in $P$, then $\{i,j\}\not\in
E$. This implies that all antichains of $P$ are stable sets of $G$, and so
$\Chn{P}\subseteq\stab{G}$. 

On the other hand, if $O$ is not transitive, then there exists two
vertices 
$k,\ell\in [n]$ such that $\{k,\ell\}\not\in E$, but such that $k$ and
$\ell$ are comparable in $P$, i.e. the transitive closure of $O$ induces
comparability of $k$ and $\ell$. Then, $e_k+e_\ell$ is a vertex of the
stable polytope $\stab{G}$ of $G$, but since $\Chn{P}$ is a subpolytope of
the $n$-dimensional cube, $e_k+e_\ell\not\in\Chn{P}$. We obtain that
$\Chn{P}\neq \stab{G}$ if $O$ is not transitive, and so $\Chn{P}\subsetneq
\stab{G}$. 

If $O$ is transitive, then $\Chn{P}=\stab{G}$. This completes the proof.
\par \qed \end{proof2}

\subsection{Poset theory.}\label{subsec:pt}

Let us now introduce the background necessary to present our second
method. This will eventually lead to a different proof of
Theorem~\ref{th:cg}.

\begin{defn}\label{def:orien}
If we consider a simple connected undirected graph $G=G(V,E)$ and endow it
with an acyclic orientation 
of its edges, we will say that our graph is an \emph{oriented graph} and consider
it a directed graph, so that every member
of $E$ is regarded as an ordered pair. We will use the notation $G_o=G_o(V,E)$ 
to denote an oriented graph defined in such a way, coming from a simple graph $G$.
\end{defn}

\begin{defn}\label{def:gf}
Let $G_o=G_o(V,E)$ be an oriented graph. We will denote by $\hat{G_o}$ the
oriented graph with vertex set $\hat{V}:=V\cup\{\hat{0},\hat{1}\}$ and set
of directed edges $\hat{E}$ equal to the union of $E$ and all edges of the
form: 
$$\begin{array}{l l l}
 & (v,\hat{1}) \text{ with $v\in V$ and $\outdeg{v}=0$ in $G_o$, and}\\
 & (\hat{0},v) \text{ with $v\in V$ and $\indeg{v}=0$ in $G_o$.} 
  \end{array}$$
A \emph{natural flow} on $G_o$ will be a function $f:\hat{E}\rightarrow\N$
such that for all $v\in V$, we have: 
$$\sum_{(x,v)\in\hat{E}}f(x,v)=\sum_{(v,y)\in\hat{E}}f(v,y).$$
In other words, a natural flow on $G_o$ is a nonnegative network flow on
$\hat{G_o}$ with unique source $\hat{0}$, unique sink $\hat{1}$, and infinite edge capacities.
\end{defn}

First, let us relate natural flows on oriented graphs with linear
extensions of induced posets.

\begin{lem}\label{lem:flow}
Let $G_o=G_o(V,E)$ be an oriented graph with induced partial order $P$ on $V$,
and with $|V|=n$.
Then, the function $g:\hat{E}\rightarrow \N$ defined by 
$$\begin{array}{l l}
g(u,v)&=\card{\set{\sigma:\sigma\text{ is a linear extension of $P$ and
$\sigma(u)=\sigma(v)-1$}}}\\
&\text{if $(u,v)\in E$,}\\
g(v,\hat{1})&=\card{\set{\sigma:\sigma \text{ is a linear extension of $P$
and $\sigma(v)=n$}}}\\
&\text{if $v\in V$ and $\outdeg{v}=0$ in $G_o$, and}\\
g(\hat{0},v)&=\card{\set{\sigma:\sigma \text{ is a linear extension of $P$
and $\sigma(v)=1$}}}\\
&\text{if $v\in V$ and $\indeg{v}=0$ in $G_o$,}
\end{array}$$
is a natural flow on $G_o$. Moreover, the net $g$-flow from 
$\hat{0}$ to $\hat{1}$ is equal to $e(P)$. 
\end{lem}
\begin{proof2}
Assume without loss of generality that $V=[n]$, and consider the directed
graph $K$ on vertex set $V(K)=[n]\cup\{\hat{0},\hat{1}\}$ whose set 
$E(K)$ of directed edges consists of all:
$$\begin{array}{r l}
(i,j) & \text{for $i<_P j$},\\
\text{$(i,j)$ and $(j,i)$} & \text{for $i||_Pj$},\\
(\hat{0},i) & \text{for $i$ minimal in $P$, and}\\
(i,\hat{1}) & \text{for $i$ maximal in $P$.}
\end{array}$$
As directed graphs, we check that $\hat{G_o}$ is a subgraph of $K$.
We will define a network flow on $K$ with unique source $\hat{0}$ and
unique sink $\hat{1}$, expressing it as a sum of simpler network flows. 

First, extend each linear extension $\sigma$ of $P$ to $V(K)$ by further
defining
$\sigma\left(\hat{0}\right)=0$ and $\sigma\left(\hat{1}\right)=n+1$. Then,
let
$f_{\sigma}:E(K)\rightarrow \N$ be given by 
$$f_{\sigma}(x,y)=\begin{cases}
1 & \text{if $\sigma(x)=\sigma(y)-1$,}\\
0 & \text{otherwise.}
\end{cases}$$
Clearly, $f_{\sigma}$ defines a network flow on $K$ with source $\hat{0}$,
sink $\hat{1}$, and total net flow $1$, and then
$f:=\displaystyle\sum_{\text{$\sigma$ linear ext. of $P$}}f_{\sigma}$
defines a network flow on $K$ with total net flow $e(P)$. Moreover, for
each $(x,y)\in \hat{E}$ we have that $f(x,y)=g(x,y)$. It remains now to
check that the restriction of $f$ to $\hat{E}$ is still a network flow on
$\hat{G_o}$ with total
flow $e(P)$. 

We have to verify two conditions. First, for $i,j\in[n]$ and if $i||_Pj$,
then $$\begin{array}{l l}
&\card{\set{\sigma:\text{ $\sigma$ is a lin. ext. of $P$ and
$\sigma(i)=\sigma(j)-1$}}}\\
=&\card{\set{\sigma:\text{ $\sigma$ is a lin. ext. of $P$ and
$\sigma(j)=\sigma(i)-1$}}},
\end{array}$$
so $f(i,j)=f(j,i)$, i.e. the net $f$-flow between $i$ and $j$ is $0$. 
Second, again for $i,j\in[n]$, if $i<_Pj$ but $i\centernot\lessdot_P j$,
then
$f(i,j)=0$. These two observations imply that $g$ defines a network flow
on $\hat{G_o}$ with total flow $e(P)$.
\par \qed \end{proof2}
The next result was obtained in~\cite{ehs} using the theory of promotion and evacuation for posets, and their proof bears
no resemblance to ours.
\begin{cor}\label{cor:tri}
Let $P$ be a partial order on $V$, with $|V|=n$. If $A$ is an antichain of
$P$, then 
$e(P)\geq \sum_{v\in A}e(P\backslash v)$, where $P\backslash v$ denotes
the
induced poset on $V\backslash\{v\}$. Similarly, if $S$ is a cutset of $P$,
then $e(P)\leq \sum_{v\in S}e(P\backslash v)$. Moreover, if $I$ is a subset
of $V$ that is either a cutset or an antichain of $P$, then $e(P)=\sum_{v\in
I}e(P\backslash v)$ if and only if $I$ is both a cutset and an antichain of
$P$.
\end{cor}
\begin{proof2}
Let $G=G(V,E)$ be any graph that contains as a subgraph the Hasse diagram
of $P$, and orient the edges of $G$ so that it induces exactly $P$ to obtain an oriented graph
$G_o$. 
Let $g$ be as in
Lemma~\ref{lem:flow}. Since edges representing cover relations of $P$ are
in $G$ and are oriented accordingly in $G_o$, the net $g$-flow is $e(P)$. Moreover, by the {\it standard chain
decomposition of network flows} of \cite{ff} (essentially
Stanley's transfer map), which expresses $g$ as a sum of positive flows
through each maximal directed path of $G_o$, it is clear that for $A$ an antichain of
$P$, we have that $e(P)\geq\sum_{v\in A}\sum_{(x,v)\in\hat{E}}g(x,v)$,
since antichains intersect maximal directed paths of $G_o$ at most once. Similarly,
for $S$ a cutset of $P$, we have that 
$e(P)\leq \sum_{v\in S}\sum_{(x,v)\in\hat{E}}g(x,v)$ since every maximal directed path 
of $G_o$ intersects $S$. Furthermore, equality will only hold in either
case if the other case holds as well. But then, for each $v\in V$, the map
$\text{{\bf Trans}}$ 
that transforms linear extensions of $P\backslash v$ into linear
extensions of $P$ and defined via: For $\sigma$ a linear extension of
$P\backslash v$ and $\kappa:=\displaystyle\max_{y<_Pv}\sigma(y)$,
$$\text{{\bf Trans}}\left(\sigma\right)(x)=\begin{cases}
\kappa+1 &\text{if $x=v$,}\\
\sigma(x)+1 &\text{if $\sigma(x)>\kappa$,}\\
\sigma(x) & \text{otherwise,}
\end{cases}$$
is a bijection onto its image, and the number
$\sum_{(x,v)\in\hat{E}}g(x,v)$ is precisely $\card{\Img{\text{{\bf
Trans}}}}$.
\par \qed \end{proof2}

Getting ready for the second proof of Theorem~\ref{th:cg}, it will be
useful to have a notation for 
the main object of study in this paper:

\begin{defn}\label{def:maxi}
Let $G=G(V,E)$ be an undirected simple graph. The maximal number of linear
extensions of a partial order on $V$ induced by an acyclic orientation of
$E$ will be denoted by $\varepsilon(G)$. 
\end{defn}
\begin{proof2}[Second proof of Theorem~\ref{th:cg}.] Assume without loss
of generality that $V=[n]$. We will do induction on $n$. The case $n=1$ is
immediate, so assume the result holds for $n-1$. Note that every induced
subgraph of $G$ is also a comparability graph and moreover, every
transitive orientation of $G$ induces a transitive orientation on the edges
of every induced graph of $G$. Now, let $O$ be a non-transitive orientation
of $E$ with induced poset $P$, so that there exists a comparable pair
$\{k,\ell\}$ in $P$ that is stable in $G$. Let $S$ be an antichain cutset
of $P$. Then, $S$ is a stable set of $G$. 
Letting $G\backslash i$ be the induced subgraph of $G$ on vertex set
$[n]\backslash\{i\}$, we obtain that
$\varepsilon(G)\geq \sum_{i\in S}\varepsilon(G\backslash i)\geq \sum_{i\in
S}
e(P\backslash i)=e(P)$,
where the first inequality is an application of Corollary~\ref{cor:tri} on a
transitive orientation of $G$, along with Definition~\ref{def:maxi} and the
inductive hypothesis, the second inequality is obtained after recognizing
that the poset induced by $O$ on each 
$G\backslash i$ is a subposet of $P\backslash i$ and by
Definition~\ref{def:maxi}, and the last equality follows because $S$ is a
cutset of $P$. If $|S|>1$ or $S\cap\{k,\ell\}=\emptyset$, then by induction
the second inequality will be strict. On the other hand, if
$S=\{k\}$ or $S=\{\ell\}$, then the first inequality will be strict since 
$\{k,\ell\}$ is stable in $G$. 

Lastly, the different posets arising from transitive orientations of $G$
have in common that their antichains are exactly the stable sets of $G$,
and their cutsets are exactly the sets that meet every maximal clique of
$G$ at least once, so by the corollary, the inductive hypothesis and our
choice of $S$ above, these posets have the same number of linear
extensions and this number is in general at least
$\sum_{i\in S}\varepsilon(G\backslash i)$, and strictly greater if
$S=\{k\}$ or $S=\{\ell\}$.   
\par \qed \end{proof2}

\section{Beyond comparability and enumerative results.}\label{sec:fe}
In this section, we will illustrate a short application of the ideas
developed in Section~\ref{sec:cg} to the case of odd cycles, 
re-establishing Theorem~\ref{th:oc} using a more elegant technique (Subsection~\ref{subsec:aut}) that applies to other families of graphs. 
Then, in Subsection~\ref{subsec:gbms},
we will obtain basic enumerative results for
$\varepsilon(G)$. Finally, in Subsection~\ref{subsec:rg}, we will 
study the random variable $\varepsilon(G)$ when $G$ is a random graph with distribution $\rgr{n}{p}$, $0<p<1$.
As it will be seen, if $G\sim\rgr{n}{p}$, then $\logt \varepsilon(G)$ concentrates tightly around its mean, and this mean 
is asymptotically equal to $n\logt\log_b n^2$, where $b=\frac{1}{1-p}$. This will permit us to
obtain, for the case of random graphs, 
new bounds for the volumes of stable polytopes, and a very strong concentration result
for the \emph{entropy} of a graph, both of which will hold $a.s.$.

\subsection{A useful technique.}\label{subsec:aut}

We start with two simple observations that remained from the theory of Section~\ref{sec:cg}.

Firstly, note that for a general graph $G$, finding $\varepsilon(G)$ is equivalent to finding the chain polytope
of maximal volume contained in $\stab{G}$, hence:
\begin{obs}\label{obs:stal}
For a simple graph $G$, we have: 
$$\varepsilon(G)\leq n!\Vol{\stab{G}}.$$
\end{obs}
Also, directly from Theorem~\ref{th:cg} we can say the following:
\begin{obs}\label{obs:help}
Let $P$ and $Q$ be partial orders on the same ground set, and suppose that
the comparability graph of $P$ contains as a subgraph the comparability
graph of $Q$. Then, $e(Q)\geq e(P)$ and moreover, if the containment of
graphs is proper, then
$e(Q)>e(P)$.
\end{obs}
\begin{proof2}[Second proof of Theorem~\ref{th:oc}.]
Note that every acyclic orientation $O$ of $E$ induces a partial order on
$V$ whose comparability graph contains (as a subgraph) the comparability
graph of a poset given by an almost bipartite orientation, and this
containment is proper if $O$ is not almost bipartite. By the symmetry of
$G$, then all of the almost bipartite orientations are equivalent.

\emph{Note to proof}: The same technique allows us to obtain results for
other restrictive families of graphs, like odd cycles with isomorphic trees
similarly attached to every element of the cycle or, perhaps more importantly, odd-anti-cycles, but we do not pursue this
here.
\par \qed \end{proof2}

\subsection{General bounds for the main statistic.}\label{subsec:gbms}

Let us now turn our attention to the general enumeration problem. Firstly, we need to dwell on the case of comparability graphs, from where
we will jump easily to general graphs. 

\begin{theo}\label{th:perf}
Let $G=G(V,E)$ be a comparability graph, and further let
$V=\{v_1,v_2,\dots,v_n\}$.
For $u_1,u_2,\dots,u_k\in V$, let $G\backslash u_1u_2\dots u_k$ be the
induced subgraph of $G$ on vertex set $V\backslash\{u_1,u_2,\dots,u_k\}$.
Then,
$$\varepsilon(G)\geq
\sum_{\sigma\in\mathfrak{S}_n}\frac{1}{\chi(G)\chi(G\backslash
v_{\sigma1})\chi(G\backslash v_{\sigma1}v_{\sigma2})\chi(G\backslash
v_{\sigma1}v_{\sigma2}v_{\sigma3})\dots \chi(v_{\sigma n})},$$
where $\mathfrak{S}_n$ denotes the symmetric group on $[n]$ and $\chi$
denotes the chromatic number of the graph. 
\end{theo}
\begin{proof2}
Let us first fix a perfect order $\omega$ of the vertices of $G$, i.g.
$\omega$
can be a linear extension of a partial order on $V$ whose comparability
graph is $G$. Let $H$ be an induced subgraph of $G$ with vertex set $V(H)$
and edge set $E(H)$, let $\omega_H$ be the restriction of $\omega$ to
$V(H)$, and let $Q$ be the partial order of $V(H)$ given by labeling every
$v\in V(H)$ with $\omega_H(v)$ and orienting $E(H)$ accordingly. Using the
colors of the optimal coloring of $H$ given by $\omega_H$, we can find
$\chi(H)$ mutually disjoint antichains of $Q$ that cover $Q$, so by
Corollary~\ref{cor:tri} we obtain that
\begin{equation}\label{eq:in}
e(Q)\geq\frac{1}{\chi(H)}\sum_{v\in V(H)}e(Q\backslash v). 
\end{equation}
Now, we note that each $Q\backslash v$ with $v\in V(H)$ is also induced by
the respective restriction of $\omega$ to $V(H)\backslash v$, 
and that the comparability of $Q\backslash v$ is $H\backslash v$, and then each of
the terms on the right hand side can be expanded similarly.
Starting from $H=G$ above and noting the fact that $\varepsilon(G)= e(Q)$
for this case, we can expand the terms of~\ref{eq:in} exhaustively to
obtain the desired expression. 
\par \qed \end{proof2}

\begin{cor}\label{cor:all}
Let $G=G(V,E)$ be any graph on $n$ vertices with chromatic number
$k:=\chi(G)$. Then $\varepsilon(G)\geq\displaystyle\frac{n!}{k^{n-k}k!}$.
\end{cor}
\begin{proof2}
We can follow the proof of Theorem~\ref{th:perf}. This time, starting from
$H=G$, $Q$ will be a poset on $V$ given by a minimal coloring of $G$, i.e.
we color $G$ using a minimal number of totally ordered colors and orient
$E$ accordingly. Then, $\varepsilon(G)\geq e(Q)$ and we can expand the
right hand side of~\ref{eq:in}, but noting that 
$Q\backslash v$ can only be guaranteed to be partitioned into at most
$\chi(G)$ antichains, and that the chromatic number of a graph is at most
the number of vertices of that graph.
\par \qed \end{proof2}

Noting that the number of cutsets is a least $2$ in most cases, a similar
argument to that of Theorem~\ref{th:perf} implies:

\begin{obs}\label{obs:cut}
Let $G=G(V,E)$ be a connected graph. Then:
$$\varepsilon(G)\leq \frac{1}{2}\sum_{v\in V}\varepsilon(G\backslash v).$$
\end{obs}
\begin{exam}
If $G=G(V,E)$ is the odd cycle on $2n+1$ vertices, then for each
$v\in V$ we have $\varepsilon(G\backslash v)=E_{2n}$, the $(2n)$-th Euler
number, and
$\chi(G)=3$, so $a_n:=\Frac{(2n+1)E_{2n}}{2}\geq \varepsilon(G)\geq
b_n:=\Frac{(2n+1)!}{3^{2n-2}\cdot3!}$. As $n$ goes to infinity, then
$\Frac{a_n}{b_n}\sim\Frac{4}{3\pi}\left(\Frac{6}{\pi}\right)^{2n}$.
\end{exam}

Other upper bounds can be obtained from rather different considerations. 
\begin{prop}\label{prop:ubb}
Let $G=G(V,E)$ be a simple graph on $n$ vertices. Then, $\varepsilon(G)$ is at most equal to the number of acyclic
orientations of the edges of $\bar{G}$, the complement of $G$. Equality is attained if and only if 
$G$ is a complete $p$-partite graph, $p\in[n]$.
\end{prop}
\begin{proof2}
Let $\bar{E}$ be the set of edges of $\bar{G}$, so that
$E\sqcup \bar{E}=\binom{V}{2}$. 

The inequality holds since two different linear extensions (understood as labelings of $V$ with
the totally ordered set $[n]$) 
of the same acyclic orientation of $E$ induce different
acyclic orientations of $\binom{V}{2}=E\sqcup\bar{E}$: As both induce the same orientation of $E$, they must induce 
different orientations of $\bar{E}$.  

To prove the equality statement, first note that if $G$ is not a complete $p$-partite graph, then there exist edges
$\{a,b\},\{a,c\}\in\bar{E}$ such that $\{b,c\}\in E$. Suppose that $(b,c)$ is a directed edge in an optimal orientation $O$ of 
$E$. Then, if we label the vertices of $\bar{G}$ with the (totally ordered) set $[n]$ 
in such a way that $c<a<b$ comparing vertices according to their 
labels, 
our labeling induces an acyclic orientation
of $\bar{E}$ which cannot be obtained from a linear extension of $O$. Hence, $\varepsilon(G)$ is strictly less than the number of
acyclic orientations of $\bar{E}$. 

If $G$ is a complete $p$-partite graph, then suppose that there exists an acyclic orientation $\bar{O}$ of $\bar{E}$ that cannot be obtained from a linear extension
of $O$, where $O$ is any optimal orientation of $E$.  Then, in the union of the (directed) edges in both $O$ and $\bar{O}$, we can find a directed cycle  
that uses at least one (directed) edge from both $O$ and $\bar{O}$. 
Take one such directed cycle with minimal number of (directed) edges. As $G$ is a comparability graph, then $O$ is transitive, and so 
the directed cycle has the form
$E_1P_1E_2P_2\dots E_mP_m$, where $E_i$ is a directed edge in $O$, $P_i$ is a directed path in $\bar{O}$, and $m\geq 1$. Let
$E_1=(a,b)$, and let $(b,c)$ be the first directed edge in $P_1$ along the directed cycle. Since $G$ is complete $p$-partite, then
$\{a,c\}\in E$ because $\{b,c\}\in\bar{E}$. Since $O$ is transitive, 
$(a,c)$ must be a directed edge in $O$. However, this contradicts the minimality of the directed
cycle. 
\par \qed \end{proof2}

\subsection{Random graphs.}\label{subsec:rg}

Changing the scope towards probabilistic models of graphs, specifically to $\rgr{n}{p}$, we will obtain a tight concentration result for the random variable
$\varepsilon(G)$ with $G\sim \rgr{n}{p}$. The central idea of the argument will be to choose an acyclic orientation of a graph $G\sim \rgr{n}{p}$ from
a minimal proper coloring of its vertices. We expect this orientation to be nearly optimal. 

Let us first recall two remarkable results that will be essential in our proof. The first one is a well-known
result of Bollob\'as, later improved on by McDiarmid:
\begin{theo}[\cite{bol},\cite{mcdiar}]\label{theo:bol}
Let $G\sim \rgr{n}{p}$ with $0<p<1$, and define $b=\Frac{1}{1-p}$. Then:
$$\chi(G)=\Frac{n}{2\log_bn-2\log_b\log_bn+O(1)}\text{ a.s.,}$$
where $\chi(G)$ is the chromatic number of $G$.
\end{theo}

To state the second result, we first need to introduce the concept of \emph{entropy} of a \emph{convex corner}, originally
defined in~\cite{lov}. We only present here the statement for the case of stable polytopes of graphs.
\begin{defn}\label{defn:entropy}
Let $G=G([n],E)$ be a simple graph, and let $\stab{B}$ be the stable polytope of $G$. Then, the \emph{entropy} 
$H(G)$ of $G$ is the quantity:
$$H(G):=\min_{a\in\stab{G}}-\sum_{i=1}^{n}\frac{1}{n}\logt a_i.$$ 
\end{defn}

In 1995, Kahn and Kim proved certain bounds for the volumes of convex corners in terms of
their entropies. One of them, when applied to 
stable polytopes, reads as follows: 
\begin{theo}[\cite{kahn}]\label{theo:kahn}
Let $G=G([n],E)$ be a simple graph, and let $\stab{G}$ be the stable polytope of $G$. Then:
$$n^n2^{-nH(G)}\geq n!\Vol{\stab{G}}\geq n!2^{-nH(G)}.$$
\end{theo}

Equipped now with these background results, the following is true:
\begin{theo}\label{theo:probg}
Let $G\sim \rgr{n}{p}$ with $0<p<1$, $b=\frac{1}{1-p}$, and write 
$s=2\log_bn-2\log_b\log_bn$. Then: 
$$\logt\varepsilon(G)\sim n\logt s\text{ holds a.s..}$$
Also, $\expe{\logt\varepsilon(G)}\sim n\logt s$. 
\end{theo}
\begin{proof2}
Let $n$ tend to infinity. Consider the chromatic number of the graph $G\sim\rgr{n}{p}$, and color $G$ properly using
$k=\chi(G)$ colors, say with color partition $a_1+a_2+\dots+a_k=n$. Then
$\logt \varepsilon(G)\geq \logt a_1!+\dots+\logt a_k!\geq k\logt  \lfloor \frac{n}{k} \rfloor!$. By Theorem~\ref{theo:bol}, we know that 
$k= \Frac{n}{s+O(1)}$ a.s., so:
\begin{equation}\label{eqn:in1}
 \logt \varepsilon(G)\geq n\logt s-\frac{n}{\ln 2}+\frac{n}{2s}(\logt s)+O\left(\frac{n}{s}\right)\text{ a.s..}
\end{equation}
We remark here that inequality~\ref{eqn:in1} gives a slightly better bound than the one obtained directly from
Corollary~\ref{cor:all}.  

Now,  the function $\logt \varepsilon$ satisfies the {\it edge Lipschitz condition} in the {\it edge exposure martingale} since
addition of a single edge to $G$ can alter $\varepsilon$ by a factor of at most $2$, so we can apply 
Azuma's inequality to obtain:
$$\proba{\card{\logt\varepsilon(G)-\expe{\logt \varepsilon(G)}}>n\left(\logt\log_b n\right)^{\frac{1}{2}}}<\Frac{2}{\log_bn}.$$ 
Combining these two results, we see that: 
$$\expe{\logt\varepsilon(G)}\geq (n\logt s)(1+o(1)),$$ 
and moreover, that
$\logt\varepsilon(G)\sim\expe{\logt \varepsilon (G)}$ a.s. holds. 

The second necessary inequality comes, firstly, from using Observation~\ref{obs:stal}, so that $\varepsilon(G)\leq n!\Vol{\stab{G}}$, and then from
a direct application of Theorem~\ref{theo:kahn}. We obtain
that $n(\logt n-H(G))\geq \logt \varepsilon(G)$. 
Now, we further observe that for $a\in \stab{G}$, we have
$\sum_i \frac{1}{n}a_i\leq \frac{1}{n}\alpha(G)$, and then:
$$H(G)=\sum_i\frac{1}{n}\left(-\logt a_i\right)\geq 
-\logt \left(\sum_i\frac{1}{n}a_i\right)\geq -\logt \frac{1}{n}\alpha(G)=\logt\frac{n}{\alpha(G)}.$$ 
A classic result of~\cite{diar} 
states that $\alpha(G)\leq s+c$ holds a.s., where
$c=2\log_b\frac{e}{2}+1$. 
Hence, a.s., $H(G)\geq \left(\logt\frac{n}{s+O(1)}\right)=\logt n-\logt(s+O(1))$, and then  
$n\logt(s+O(1))\geq \logt\varepsilon(G)$. From here, we directly obtain:
\begin{equation}\label{eqn:in2}
\logt\varepsilon(G)\leq n\logt s+O\left(\frac{n}{s}\right)\text{ a.s.}.
\end{equation}

Therefore, from inequalities~\ref{eqn:in1} and~\ref{eqn:in2}:
$$\logt\varepsilon(G)= n\logt s+O(n)\text{ a.s..}$$  
\par \qed \end{proof2}
Calculating inequality~\ref{eqn:in2} more precisely by dropping the $O$-notation and using
Grimmett and McDiarmid's constant, we obtain:
\begin{cor}\label{cor:stab}
Let $G\sim \rgr{n}{p}$ with $0<p<1$, $b=\frac{1}{1-p}$ and $s=2\log_bn-2\log_b\log_bn$. Then, for large enough $n$:
$$\Frac{s^n}{n!}\cdot \left(\frac{1}{e}\right)^n\leq\Vol{\stab{G}}\leq \Frac{s^n}{n!}\cdot c^{n/s}\text{ a.s., where $c=2\left(\Frac{e}{2}\right)^{2/(\logt b)}$}.$$
\end{cor}
\begin{cor}\label{cor:entropy}
Let $G\sim \rgr{n}{p}$ with $0<p<1$, $b=\frac{1}{1-p}$ and $s=2\log_bn-2\log_b\log_bn$. Then, for large enough $n$:
$$\logt\left(\frac{n}{s}\right)+O\left(\frac{1}{s}\right)  \leq H(G)\leq \logt\left(\frac{n}{s}\right)+\frac{1}{\ln 2}\text{ a.s.}.$$
\end{cor}

\section{Further techniques.}\label{sec:ft}
In this section, we will see how the main problem has two more presentations as
selecting a region in the graphical arrangement with maximal
\emph{fractional volume}, or as selecting a vertex of the graphical
zonotope that is farthest from the origin in Euclidean distance.  

\begin{defn}\label{def:ga}
Consider a simple undirected graph $G=G([n],E)$. The \emph{graphical arrangement}
of $G$ is the central hyperplane arrangement in $\R^n$ given by:
$$\Arr_G=\{x\in\R^n:x_i-x_j=0\text{ , $\forall$ $\{i,j\}\in E$}\}.$$
\end{defn}
   
The regions of the graphical arrangement $\Arr_G$ with $G=G([n],E)$ are in
one-to-one correspondence with the acyclic orientations of $G$. Moreover, the
complete fan in $\R^n$ given by $\Arr_G$ is combinatorially dual to the
{\it graphical zonotope} of $G$: 
$$\mathcal{Z}_{G}^{\scriptscriptstyle\emph{central}} := \displaystyle\sum_{\{i,j\}\in
E}\left[e_i-e_j,e_j-e_i\right],$$
and there is a clear correspondence between the regions of $\Arr_{G}$ and
the vertices of $\mathcal{Z}_{G}^{\scriptscriptstyle\emph{central}}$.

Following \cite{ks}, we define the \emph{fractional volume} of a region
$\mathcal{R}$ of $\Arr_G$ to be:
$\Voli{\mathcal{R}}=\Frac{\Vol{B^n\cap \mathcal{R}}}{\Vol{B^n}}$,
where $B^n$ is the unit $n$-dimensional ball in $\R^n$. 

With little work it is possible to say the following about these volumes:

\begin{prop}\label{prop:vols}
Let $G=G([n],E)$ be an undirected simple graph, and let $\Arr_G$ be its
graphical arrangement. If $\mathcal{R}$ is a region of $\Arr_G$ and $P$ is
its corresponding partial order on $[n]$, then:
$$\Voli{\mathcal{R}}=\frac{e(P)}{n!}.$$
\end{prop}

The problem of finding the regions of $\Arr_G$ with maximal fractional volume is, intuitively, 
closely related
to the problem of finding the vertices of $\mathcal{Z}_{G}^{\scriptscriptstyle\emph{central}}$ that are farthest from the origin under some
appropriate choice of metric. It turns out that, with Euclidean metric, a precise statement can be formulated when $G$ is a comparability
graph:

\begin{theo}\label{theo:cgz}
Let $G=G(V,E)$ be a comparability graph. Then, the vertices of the
graphical zonotope of $\mathcal{Z}_{G}^{\scriptscriptstyle\emph{central}}$ that have maximal 
Euclidean distance to the origin are precisely those that correspond to the
transitive orientations of $E$, which in turn have maximal number
$\varepsilon(G)$ of linear extensions.  
\end{theo} 

To prove Theorem~\ref{theo:cgz}, we first 
note that for a simple (undirected) graph $G=G(V,E)$, the vertex of $\mathcal{Z}_{G}^{\scriptscriptstyle\emph{central}}$ corresponding 
to a given acyclic orientation of $E$ is precisely the point:
$$(\outdeg{v}-\indeg{v})_{v\in V},$$
where $\outdeg{\cdot}$ and $\indeg{\cdot}$ are calculated using the given orientation.

We need to establish a preliminary lemma. 
\begin{lem}\label{lem:inout}
Let $G_o=G_o(V,E)$ be an oriented graph. Then,
$$\frac{1}{2}\sum_{v\in
V}\left(\indeg{v}-\outdeg{v}\right)^2=|E|+\tri{G_o}+\inc{G_o}-\com{G_o},$$
where: 
\begin{itemize}
\item[1.] $\tri{G_o}$ is the number of directed triangles $(u,v),(v,w),(u,w)\in
E$. 
\item[2.] $\inc{G_o}$
is the number of triples $u,v,w\in V$ such that $(v,w),(w,v)\not\in
E$ but either $(u,v),(u,w)\in E$ or $(v,u),(w,u)\in E$. 
\item[3.] $\com{G_o}$ is
the number of
directed $2$-paths $(u,v),(v,w)\in E$ such that $(u,w)\centernot\in E$. 
\end{itemize}
\end{lem}
\begin{proof2}
For $v\in V$, $\outdeg{v}^2$ is equal to $\outdeg{v}$ plus two times the
number of pairs $u\neq w$ such that $(v,u),(v,w)\in E$, $\indeg{v}^2$ is
equal to $\indeg{v}$ plus two times the number of pairs $u,\neq w$ such
that $(u,v),(w,v)\in E$, and $\outdeg{v}\cdot\indeg{v}$ is equal to the number
of pairs 
$u\neq w$ such that $(u,v),(v,w)\in E$. If we add up these terms and
cancel out terms in the case of directed triangles, we obtain the desired
equality.   
\par \qed \end{proof2}
An important consequence of Lemma~\ref{lem:inout} is the following: 
\begin{itemize}
\item[] If $G=G(V,E)$ is a simple graph, all the
acyclic orientations of $E$ will not vary in their values of $\tri{\cdot}$
and of
$|E|$, which depend on $G$, 
but only in $\com{\cdot}$ and $\inc{\cdot}$. Moreover, $\com{\cdot}+\inc{\cdot}$ is 
equal to the number of $2$-paths in $G$ of the form
$\{u,v\},\{v,w\}\in E$ with $u\neq w$, so it is also independent
of the choice of orientation for $E$.
\end{itemize}
\begin{proof2}[Proof of Theorem~\ref{theo:cgz}.] 
We apply Lemma~\ref{lem:inout} directly.
Since $G$ is a comparability graph, from Theorem~\ref{th:cg}, 
we know that the value of $\inc{\cdot}-\com{\cdot}$ will be maximized precisely on the
transitive orientations of $G$, since all transitive orientations force $\com{\cdot}=0$. 
\par \qed \end{proof2}

\begin{proof3}[Acknowledgements:]
I am thankful to my advisor Richard P. Stanley for his support and encouragement, and to Carly Klivans and 
Federico Ardila for useful comments during the edition stage of this work. Thanks are owed to two anonymous
reviewers at FPSAC for further edition suggestions, and to the Research Science Institute of MIT for monetary 
support.  
\end{proof3}

\nocite{*}
\bibliographystyle{abbrvnat}
\bibliography{tt}
\label{sec:biblio}

\end{document}